\documentclass[11pt]{article}
\usepackage[a4paper,left=2.8cm,right=2.8cm,top=2.5cm,bottom=2.5cm]{geometry}
\usepackage{amsmath,amssymb,amsthm,mathtools,bm,booktabs,array,graphicx}
\usepackage{placeins}
\usepackage{xcolor}
\definecolor{Mycolor}{RGB}{10,74,38}
\usepackage[colorlinks=true,linkcolor=Mycolor,citecolor=Mycolor,urlcolor=Mycolor]{hyperref}
\newtheorem{theorem}{Theorem}[section]

\newtheorem{assumption}{Assumption}[section]
\newtheorem{lemma}{Lemma}[section]

\theoremstyle{remark}
\newtheorem{remark}{Remark}[section]
\numberwithin{equation}{section}
\newcommand{\Sph}{\mathbb S^2}
\newcommand{\R}{\mathbb R}
\newcommand{\dd}{\,\mathrm d}
\newcommand{\Hs}{H^s(\Sph)}

\title{Graded cubature with quasi-Monte Carlo designs for weakly singular integrals on the sphere}
\author{Hao-Ning Wu\thanks{School of Mathematical Sciences and Fujian Key
Laboratory of Mathematical Modeling and High Performance Scientific Computing,
Xiamen University, Xiamen, Fujian 361005, China
(\texttt{hnwu@xmu.edu.cn}).}}
\date{}

\begin{document}
\maketitle
\begin{abstract}
Weakly singular integrals over smooth surfaces, when pulled back to the sphere, generally have direction-dependent singularities. We study their numerical evaluation by combining radial grading transformations with spherical quasi-Monte Carlo (QMC) designs. For kernels of singularity order $0<\alpha<2$ satisfying scale-dependent derivative bounds, we prove that a grading exponent $q$ satisfying $q(2-\alpha)>s+1$ preserves the $H^s(\Sph)$ regularity of the density in the transformed integrand, uniformly in the target, for every real $s>1$. The condition is shown to be sharp for the kernel class considered. Consequently, $N$-point QMC designs for $H^s(\Sph)$ retain their $O(N^{-s/2})$ cubature error rate for the graded integrals. In addition to covering two classical transformations of Atkinson and Sidi, our analysis accommodates a localized cap-area transformation proposed here, which has a simple surface Jacobian and leaves the sphere unchanged outside a prescribed cap. Numerical experiments illustrate the theoretical results
and the applicability of graded cubature to boundary and close evaluation of Laplace double-layer potentials on a nonconvex surface.
\end{abstract}


\section{Introduction}\label{sec:intro}
Integrability and suitability for high-order numerical integration are two rather different properties of a function. Weakly singular integrals encountered in potential theory and boundary integral equations provide an illustration of this distinction. Let $\Gamma\subset\mathbb R^3$ be a smooth closed surface admitting a smooth diffeomorphism $\Phi:\Sph\to\Gamma$, where $\Sph$ denotes the unit sphere. Consider a weakly singular integral operator on $\Gamma$ of the form
\[
 (\mathcal A\mu)(X)
 =\int_\Gamma
 \frac{a(X,Y)}{|X-Y|^\alpha}\mu(Y)\dd S_\Gamma(Y),
 \qquad 0<\alpha<2,
\]
where $a$ is smooth. Denoting by $X=\Phi(x)$, $Y=\Phi(y)$, and $f=\mu\circ\Phi$, we pull the integral back to the sphere and obtain
\begin{equation}\label{eq:intro-model}
 (\mathcal Kf)(x)=\int_{\Sph}K(x,y)f(y)\dd S(y),
\end{equation}
where
\[
 K(x,y)=
 \frac{a(\Phi(x),\Phi(y))j_\Phi(y)}
 {|\Phi(x)-\Phi(y)|^\alpha},
\]
$j_\Phi$ is the surface Jacobian of $\Phi$, and $\dd S$ denotes surface measure on $\Sph$, with $|\Sph|=4\pi$. Pulling an integral back to the sphere simplifies the integration domain at the expense of complicating the kernel. In particular, even for a smooth amplitude $a$, the singularity of $K$ need not be carried by a zonal factor with a smooth remaining amplitude. Its leading coefficient may depend on the direction of approach to the diagonal. Although the integral is finite for continuous densities $f$, its integrand may be unbounded at $y=x$, even when $f$ is infinitely differentiable.

Standard cubature estimates for smooth integrands are therefore not directly applicable. A related difficulty arises when the target lies close to, rather than on, the surface: the integrand may develop a narrow peak that becomes increasingly difficult to resolve as the target approaches the boundary. Various approaches to singular quadrature have been developed, differing mainly in how the singularity is treated. Rotation of the singularity to a pole facilitates product integration in polar coordinates \cite{Atkinson1982Sphere,GrahamSloan2002}. For kernels whose singularities are carried by zonal factors, spherical product integration also provides an effective means of discretizing second-kind integral equations \cite{MR5025409}. An alternative is to modify the integration coordinates. Radial grading and angular endpoint transformations weaken the singularity through a vanishing Jacobian, such as the Atkinson transform \cite{Atkinson2004Surface} and Sidi transform \cite{Sidi2005Sphere,Sidi2006Surface}. The convergence of the resulting product quadrature rules has been analyzed through endpoint expansions \cite{AtkinsonSommariva2005,Sidi2005Atkinson}. On local patches, the Duffy transformation offers a related mechanism for absorbing singular factors into the Jacobian \cite{Duffy1982}. Other approaches act more directly on the integrand or the quadrature rule: singularity subtraction and density interpolation remove local singular contributions \cite{FariaPerezBonnet2021},
corrected quadrature modifies weights near the singularity \cite{WuMartinsson2023}, and quadrature by expansion uses local expansions of layer potentials away from the boundary \cite{KlocknerEtAl2013}. Singularity reconstruction improves finite element convergence for fractional boundary value problems by separating a singular component from a more regular remainder \cite{JinZhou2015}. We focus on changes of variables that improve the regularity of the integrand before applying a cubature rule. The question is how much regularity grading provides, independently of the cubature rule.

Once a singularity has been sufficiently weakened by a change of variables, the choice of a cubature rule becomes largely a question of how well the underlying point set integrates regular functions. On the sphere, a classical notion of integration quality is \emph{polynomial exactness}. A spherical $t$-design \cite{delsarte1991geometriae}  integrates every spherical polynomial of degree at most $t$ exactly with equal weights. The existence of spherical $t$-designs with $N\asymp t^2$ points was established in \cite{bondarenko2013optimal}. For functions in the intrinsic Sobolev space $H^s(\Sph)$ with $s>1$, such designs achieve an integration error of order $N^{-s/2}$; see, e.g., \cite{BrauchartHesse2007,MR3246811,HesseSloan2005Worst,HesseSloan2006Sobolev}.

Polynomial exactness, however, is not a prerequisite for achieving this rate. For a fixed $s>1$, a sequence of point sets is called a sequence of quasi-Monte Carlo (QMC) designs \cite{MR3246811} for $H^s(\Sph)$ if the associated equal-weight cubature rules
\begin{equation}\label{equ:quad}
 Q_Nf=\frac{4\pi}{N}\sum_{i=1}^Nf(z_{N,i})
\end{equation}
satisfy
\begin{equation}\label{eq:intro-qmc}
 \text{wce}(Q_N;H^s)
 :=\sup_{\|f\|_{H^s(\Sph)}\le1}
 \left|\int_{\Sph}f\dd S-Q_Nf\right|
 \le C_QN^{-s/2},
\end{equation}
where $C_Q$ is independent of $N$.
Consequently, for every $f\in H^s(\Sph)$,
\[
 \left|\int_{\Sph}f\dd S-Q_Nf\right|
 \le \text{wce}(Q_N;H^s)\|f\|_{H^s(\Sph)}
 \le C_QN^{-s/2}\|f\|_{H^s(\Sph)}.
\]
The restriction $s>1$ ensures the Sobolev embedding $H^s(\Sph)\hookrightarrow C(\Sph)$, thus point evaluations in the cubature rule \eqref{equ:quad} are well defined. Unlike spherical $t$-designs, QMC designs are characterized by their
worst-case integration errors in Sobolev spaces rather than by exactness on polynomial spaces; see also \cite{SloanWomersley2026QMC,WuZhuang2026NestedQMC,Wu2026PerturbedQMC}. For PDEs with random inputs, weighted Sobolev regularity estimates similarly underpin QMC convergence analysis \cite{WuZhangZhao2026}. We therefore study the Sobolev regularity of the transformed integrand.

To apply this Sobolev error estimate to weakly singular integrals, the transformed integrand must have sufficient regularity. Changes of variables combined with spherical $t$-designs have been explored numerically
\cite{AnChen2016}, without a corresponding convergence analysis. The question we address here is how the Sobolev regularity of the transformed integrand depends on the singularity order of the kernel, the grading exponent, and the regularity of the density. This requires an analysis of the transformed
integrand as a whole, rather than merely the cancellation of the kernel singularity by the Jacobian. In particular, we seek estimates uniform in the target point for densities of finite Sobolev regularity.

The idea is to concentrate the cubature points near the singularity without having to construct a new point set for each target. This can be achieved by contracting geodesic distances towards the singularity and compensating for the contraction through the surface Jacobian. The question is how much contraction is needed to make the Sobolev error estimate \eqref{eq:intro-qmc} applicable. To illustrate the idea, we fix a target $x\in\Sph$ and write
a point $z$ near $x$ in geodesic polar coordinates, with $r=d(x,z)$ denoting the geodesic distance on $\Sph$, as
\[
 z=\cos r\,x+\sin r\,\omega,
\]
where $\omega\in T_x\Sph$ satisfies $|\omega|=1$. For a grading exponent $q>1$, we define a radial transformation $\mathcal{T}_{x,q}$ by
\[
 \mathcal{T}_{x,q}z
 =\cos\tau_q(r)\,x+\sin\tau_q(r)\,\omega,
\]
where $\tau_q:[0,\pi]\to[0,\pi]$ is a strictly increasing
radial profile that agrees with $r^q$ near $r=0$ and
with the identity outside a neighbourhood of the origin.
Thus the transformation leaves the angular direction unchanged
but contracts radial distances near $x$.
With $j_\Phi$ already included in $K$, let $J_{x,q}$ denote the surface
Jacobian of $\mathcal{T}_{x,q}$ with respect to $\dd S$.
Changing variables in \eqref{eq:intro-model} gives
\begin{equation}\label{eq:intro-g}
 (\mathcal Kf)(x)
 =\int_{\Sph}G_{x,q}f(z)\dd S(z),
 \end{equation}
 where 
 \[ G_{x,q}f(z):=K(x,\mathcal{T}_{x,q}z)f(\mathcal{T}_{x,q}z)J_{x,q}(z).\]
The local behaviour of the Jacobian explains how grading weakens the singularity:
\[
 J_{x,q}(z)
 =\frac{\tau_q'(r)\sin\tau_q(r)}{\sin r}
 \sim q r^{2q-2}
 \qquad\text{as }r\to0.
\]
Near $r=0$, the transformed distance is
$d(x,\mathcal{T}_{x,q}z)=\tau_q(r)=r^q$.
Since $\Phi$ is a smooth regular parametrization,
$|\Phi(x)-\Phi(\mathcal{T}_{x,q}z)|\asymp r^q$ locally.
The factors $a$ and $j_\Phi$ in the kernel are bounded, so
\[
 K(x,\mathcal{T}_{x,q}z)=O(r^{-q\alpha}).
\]
Combining this with the Jacobian asymptotic gives
\[
 K(x,\mathcal{T}_{x,q}z)J_{x,q}(z)
 =O\bigl(r^{q(2-\alpha)-2}\bigr).
\]
A sufficiently large $q$ therefore makes the transformed
weight vanish at the singularity. Continuity alone, however,
does not suffice for high-order cubature. What matters for
\eqref{eq:intro-qmc} is the Sobolev regularity of the entire
transformed integrand, including the density composed with
the grading transformation.

Our main result identifies a sharp threshold for the grading exponent. For any real $s>1$, the grading condition
\begin{equation}\label{eq:intro-grading}
 q(2-\alpha)>s+1
\end{equation}
ensures that
\begin{equation}\label{eq:intro-bound}
 \sup_{x\in\Sph}\|G_{x,q}f\|_{H^s(\Sph)}
 \le C\|f\|_{H^s(\Sph)},
 \qquad f\in H^s(\Sph),
\end{equation}
where $C$ may depend on the kernel and the grading profile, but is independent of the target $x$ and the density $f$. In particular, the transformed integrand inherits the Sobolev regularity of the density without requiring any additional smoothness of $f$. We will prove this estimate \eqref{eq:intro-bound} for kernels whose source derivatives of order $j$ are bounded by $C d(x,y)^{-\alpha-j}$, uniformly in $x$, for $0\le j\le\lfloor s\rfloor+1$. The strict inequality in \eqref{eq:intro-grading} cannot be relaxed for this kernel class. When $q(2-\alpha)\le s+1$, we construct a positive kernel with a direction-dependent singularity for which the transformed integrand fails to belong to $H^s(\Sph)$, even when the density is constant.

The significance of \eqref{eq:intro-bound} for numerical integration is that the analysis of the singularity can be separated from that of the cubature rule. Given a QMC design for $H^s(\Sph)$, we define the \emph{graded cubature rule}\footnote{The term ``graded cubature'' is motivated by the
terminology of Atkinson \cite{Atkinson2004Surface} and Sidi \cite{Sidi2005Atkinson}, who refer to $q$ as the ``grading parameter'' in their studies of variable transformations for numerical integration over surfaces.} by
\begin{equation}\label{eq:graded-rule}
 (\mathcal K_N^{\rm gr}f)(x)
 :=Q_N(G_{x,q}f).
\end{equation}
Under the grading condition \eqref{eq:intro-grading}, the transformed integrand admits a continuous extension with value zero at the target $x$. Combining \eqref{eq:intro-qmc} and \eqref{eq:intro-bound} gives
\begin{equation*}
 \|\mathcal Kf-\mathcal K_N^{\rm gr}f\|_\infty\le C_QN^{-s/2}
 \sup_{x\in\Sph}\|G_{x,q}f\|_{H^s(\Sph)}\le C N^{-s/2}\|f\|_{H^s(\Sph)}.
\end{equation*}
Thus grading gives the required Sobolev regularity, and the QMC-design bound gives the integration error.

The rest of this paper is organized as follows. Section~\ref{sec:preliminaries} reviews Sobolev spaces and states the kernel assumptions. Section~\ref{sec:grading} introduces grading transformations, including the localized cap-area construction. Section~\ref{sec:regularity} establishes the sharp grading threshold for Sobolev regularity. Section~\ref{sec:cubature} deduces cubature error bounds from the regularity theorem. Section~\ref{sec:numerics} presents numerical experiments.

\section{Preliminaries}\label{sec:preliminaries}\label{sec:setting}

We review Sobolev spaces and state the kernel assumptions.

\subsection{Sobolev spaces}\label{sec:sobolev}
Let $d(x,y):=\arccos(x\cdot y)$ be the geodesic distance between $x,y\in\Sph$, and let $\dd S$ denote surface measure on $\Sph$, normalized so that $|\Sph|=4\pi$. For $\ell\ge0$, let $\{Y_{\ell,m}:-\ell\le m\le\ell\}$ be a real orthonormal basis of spherical harmonics of degree $\ell$ in $L^2(\Sph,\dd S)$. For $f\in L^2(\Sph,\dd S)$, its spherical harmonic coefficients are given by
\[
 \widehat f_{\ell,m}
 :=\int_{\Sph}f(y)Y_{\ell,m}(y)\dd S(y).
\]
For $s\ge0$, the Sobolev space $H^s(\Sph)$ consists of all functions $f\in L^2(\Sph)$ for which
\begin{equation}\label{eq:Hs}
 \|f\|_{\Hs}^2 := \sum_{\ell=0}^\infty\sum_{m=-\ell}^{\ell} (1+\ell(\ell+1))^s|\widehat f_{\ell,m}|^2 <\infty.
\end{equation}
This norm is equivalent to the local Euclidean Sobolev norms in smooth
charts \cite{MR2934227,DaiXu2013}.
For an integer $m\ge0$, $C^m$ denotes the space of functions with continuous
derivatives through order $m$. On a compact domain $E$, we use the norm
\begin{equation*}
 \|f\|_{C^m(E)}:=\max_{0\le j\le m}\sup_{y\in E}|\nabla^j f(y)|,
\end{equation*}
where $\nabla^0f=f$, and $\nabla^j f$ denotes the $j$th covariant derivative
on $\Sph$ or the $j$th ordinary derivative in Euclidean coordinates,
with its usual pointwise tensor norm.

Throughout the paper, we assume $s>1$ and denote $k:=\lfloor s\rfloor$ and $M:=k+1$.
By the Sobolev embedding $H^s(\Sph)\hookrightarrow C(\Sph)$,
every $f\in H^s(\Sph)$ has a continuous representative
satisfying
\begin{equation}\label{eq:embedding}
 \|f\|_\infty\le C_s\|f\|_{\Hs},
\end{equation}
where $\|\cdot\|_\infty$ denotes the supremum norm on $\Sph$
and $C_s$ depends only on $s$.

\subsection{Kernel assumptions}\label{sec:kernels}
To describe the local behaviour of a weakly singular kernel, we use geodesic normal coordinates centered at the target point $x$. For each $x\in\Sph$, let $T_x\Sph=\{v\in\R^3:x\cdot v=0\}$
be the tangent space at $x$. The Riemannian exponential map $\exp_x:T_x\Sph\to\Sph$ is given explicitly by
\[
 \exp_x(v)=\cos|v|\,x+\frac{\sin|v|}{|v|}v, \qquad v\ne0,
\]
with $\exp_x(0)=x$, where $|v|$ denotes the Euclidean norm. In particular, $d(x,\exp_x(v))=|v|$ for $|v|<\pi$. We identify $T_x\Sph$ with $\R^2$ through an arbitrary orthonormal frame and use the same notation $\exp_x$ under this identification.

To control the Sobolev norm after grading, we need estimates for the singular kernel $K$ and its source derivatives on annuli shrinking to the singular point. For a kernel with a singularity of order $\alpha$, its size on an annulus of radius $t$ is $O(t^{-\alpha})$. We rescale it to a fixed annulus and multiply by $t^\alpha$ to remove this growth. We therefore require uniform $C^M$ bounds for the rescaled kernel, together with ordinary $C^M$ bounds away from the diagonal. To this end, we fix a radius $0<r_*<\pi/8$ and introduce the reference annulus
\[
 A_*:=\{v\in\R^2:1/4<|v|<4\}.
\]
Let $0<\alpha<2$ and let $K$ be a kernel that is $C^M$ in its second variable away from the diagonal $\{(x,x):x\in\Sph\}$. We measure the size of $K$ and its source derivatives by
\begin{equation}\label{eq:kernel-size}
 \|K\|_{\alpha,M}
 :=
 \sup_{x\in\Sph}
 \|K(x,\cdot)\|_{C^M(\{y\in\Sph:d(x,y)\ge \tfrac{r_*}{2}\})}+\sup_{\substack{x\in\Sph\\0<t\le r_*}}
 \|t^\alpha K(x,\exp_x(t\,\cdot))\|_{C^M(\overline{A_*})}.
\end{equation}
Here the first $C^M$ norm is defined in terms of covariant derivatives on $\Sph$, while the second is the usual Euclidean $C^M$ norm with respect to the rescaled variable $v\in\overline{A_*}$, where $\overline{A_*}=\{v\in\R^2:1/4\le |v|\le4\}$ is the closure of $A_*$. The derivatives are understood on
open neighbourhoods of the indicated closed sets, and the bounds are required to be uniform in both the target $x$ and the choice of orthonormal frame.

\begin{assumption}\label{ass:kernel}
The kernel satisfies $\|K\|_{\alpha,M}<\infty$ for some
$0<\alpha<2$.
\end{assumption}

The second term in \eqref{eq:kernel-size} is finite if and only if
there is a constant $C$, independent of $x$ and the orthonormal frame,
such that
\begin{equation}\label{eq:kernel-derivatives}
 |D_\zeta^\nu[K(x,\exp_x\zeta)]|\le C|\zeta|^{-\alpha-|\nu|},\qquad
 |\nu|\le M,\quad 0<|\zeta|<4r_*.
\end{equation}
Thus Assumption~\ref{ass:kernel} is equivalent to requiring both
finiteness of the first term in \eqref{eq:kernel-size} and the local
derivative bounds \eqref{eq:kernel-derivatives}.
Indeed, for fixed $t$, the chain rule gives
\[
 D_v^\nu\!\left[t^\alpha K(x,\exp_x(tv))\right]
 =t^{\alpha+|\nu|}
 \left.D_\zeta^\nu[K(x,\exp_x\zeta)]\right|_{\zeta=tv}.
\]
If the second term in \eqref{eq:kernel-size} is finite, for $0<|\zeta|<4r_*$ we choose $t=\min\{|\zeta|,r_*\}$ and $v=\zeta/t$. Then $1\le|v|<4$ and $t\asymp|\zeta|$, so this identity gives \eqref{eq:kernel-derivatives}. Conversely, for $v\in\overline{A_*}$, $0<t\le r_*$, and $|\nu|\le M$, \eqref{eq:kernel-derivatives} gives, by continuity at the boundary,
\[
 \left|D_v^\nu[t^\alpha K(x,\exp_x(tv))]\right|
 \le C|v|^{-\alpha-|\nu|}\le C4^{\alpha+|\nu|}.
\]
This bounds the second term in \eqref{eq:kernel-size} and proves the equivalence. The same calculation gives uniform rescaled $C^M$ bounds on any fixed annulus away from the origin, at sufficiently small scales. Taking $|\nu|=0$ in
\eqref{eq:kernel-derivatives} and using the first term in
\eqref{eq:kernel-size} away from the diagonal gives
\begin{equation}\label{eq:kernel-pointwise}
 |K(x,y)|\le C\|K\|_{\alpha,M}d(x,y)^{-\alpha},\qquad x\ne y.
\end{equation}
For continuous $f$, the integral in \eqref{eq:intro-model} is absolutely
convergent, since $\int_0^{r_*}r^{-\alpha}\sin r\dd r<\infty$ for $\alpha<2$.

Suppose the first term in \eqref{eq:kernel-size} is finite. A class of kernels satisfying Assumption~\ref{ass:kernel} is given locally in geodesic polar coordinates by
\begin{equation}\label{eq:directional}
 K(x,\exp_x(r\omega))=r^{-\alpha}b_x(r,\omega),
\end{equation}
where $b_x$ and its radial and angular derivatives through order $M$ are uniformly bounded for $0\le r\le4r_*$. Indeed, on the reference annulus $\overline{A_*}$, letting $r=t|v|$ and $\omega=v/|v|$ gives
\[
 t^\alpha K(x,\exp_x(tv))
 =|v|^{-\alpha}
 b_x\left(t|v|,\frac{v}{|v|}\right),
\]
which has uniform $C^M$ bounds.
In particular, $b_x(0,\omega)$ may depend on the direction
$\omega$, so the kernel need not admit a decomposition into
a zonal singular factor and a smooth amplitude.
For the surface pullback in \eqref{eq:intro-model}, Taylor expansion gives
\[
 b_x(0,\omega)
 =\frac{a(\Phi(x),\Phi(x))j_\Phi(x)}{|D\Phi_x\omega|^\alpha},
\]
where $D\Phi_x$ is the differential of $\Phi$ at $x$.
Its dependence on $\omega$ comes from the directional stretching of
the parametrization.

\section{Grading transformations}\label{sec:grading}
For an arbitrary target $x\in\Sph$, let $z=\cos\theta\,x+\sin\theta\,\omega$ for $0<\theta<\pi$, where $\theta=d(x,z)$ and $\omega\in T_x\Sph$ with $|\omega|=1$. Let $\tau:[0,\pi]\to[0,\pi]$ be a continuous function satisfying $\tau(0)=0$ and $\tau(\pi)=\pi$. We assume that $\tau$ is smooth on $(0,\pi)$ with $\tau'>0$. The radial transformation associated with $\tau$ is defined by
\begin{equation}\label{eq:radial-map}
 \mathcal{T}_{x,\tau}(z)
 =\cos\tau(\theta)\,x+\sin\tau(\theta)\,\omega.
\end{equation}
By continuity, $\mathcal{T}_{x,\tau}$ extends to a
homeomorphism of $\Sph$ onto itself that fixes both $x$
and $-x$. Its surface Jacobian is given by
\begin{equation}\label{eq:surface-J}
 J_{x,\tau}(z)
 =\tau'(\theta)\frac{\sin\tau(\theta)}{\sin\theta},
 \qquad 0<\theta<\pi.
\end{equation}
Since $\dd S(z)=\sin\theta\,\dd\theta\,\dd\phi$ in spherical
coordinates, where $\phi$ denotes the azimuthal angle, we have
\[
 J_{x,\tau}(z)\dd S(z)
 =\tau'(\theta)\sin\tau(\theta)\,\dd\theta\,\dd\phi.
\]

\begin{assumption}\label{ass:grading}
Let $\tau$ be a radial profile satisfying the conditions stated above, and let $x\in\Sph$ be an arbitrary target. We assume that there exist $q>1$, $0<r_0<\pi/2$, and positive constants $c$, $C$, and $C_j$ such that $\tau$ satisfies
\begin{equation}\label{eq:profile-bounds}
 \begin{aligned}
 c r^q &\le \tau(r)\le C r^q,\\
 c r^{q-1} &\le \tau'(r)\le C r^{q-1},\\
 |\tau^{(j)}(r)| &\le C_j r^{q-j},
 \qquad 2\le j\le M+1,
 \end{aligned}
\end{equation}
for all $0<r<r_0$, where $r=d(x,z)$ is the distance to the target. At its antipode $-x$, now letting $r=d(-x,z)$, we assume that either $\tau(\pi-r)=\pi-r$ for $0<r<r_0$, so that the transformation is the identity in a neighbourhood of $-x$, or the profile $\tau_-(r):=\pi-\tau(\pi-r)$
satisfies the bounds in \eqref{eq:profile-bounds} with the same exponent $q$, possibly with different constants. The radius $r_0$ and all constants in the above bounds are independent of the target $x$.
\end{assumption}

The first two bounds in \eqref{eq:profile-bounds} prescribe the strength of the grading, while the remaining bounds control its higher-order derivatives. The assumption permits both transformations localized near the singularity and transformations that grade distances towards both the target and its antipode.

\subsection{Classical grading transformations}\label{sec:classical-grading}

We first recall two classical transformations for the numerical evaluation of singular integrals.

\textbf{Atkinson's transformation.} For $q>1$ and the target at the standard north pole $(0,0,1)$, Atkinson's transformation, proposed in \cite{Atkinson2004Surface}, is given by
\begin{equation}\label{eq:Atkinson}
 (\sin\theta\cos\phi,\sin\theta\sin\phi,\cos\theta)
 \longmapsto
 \frac{(\sin^q\theta\cos\phi,\sin^q\theta\sin\phi,\cos\theta)}
 {\sqrt{\cos^2\theta+\sin^{2q}\theta}}.
\end{equation}
The radial profile $\tau_A(\theta)\in[0,\pi]$ is determined by
\[
 \cos\tau_A(\theta)
 =\frac{\cos\theta}{\sqrt{\cos^2\theta+\sin^{2q}\theta}}.
\]
Since $\tau_A(\theta)\in[0,\pi]$, we also have
\[
 \sin\tau_A(\theta)
 =\frac{\sin^q\theta}
 {\sqrt{\cos^2\theta+\sin^{2q}\theta}}.
\]
Substitution into \eqref{eq:radial-map}, with $x=(0,0,1)$, recovers \eqref{eq:Atkinson}. The derivative and surface Jacobian of $\tau_A(\theta)$ satisfy
\begin{equation}\label{eq:Atkinson-derivative}
 \tau_A'(\theta)
 =\frac{\sin^{q-1}\theta(q\cos^2\theta+\sin^2\theta)}
 {\cos^2\theta+\sin^{2q}\theta}
\end{equation}
and
\begin{equation*} 
 J_A(\theta)
 =\frac{\sin^{2q-2}\theta(q\cos^2\theta+\sin^2\theta)}
 {(\cos^2\theta+\sin^{2q}\theta)^{3/2}}.
\end{equation*}
Near $r=0$, we have
\[
 \tau_A(r)=\arctan\left(\frac{\sin^q r}{\cos r}\right).
\]
Together with $z=\sin^q r/\cos r$, we obtain
\[
 \frac{\tau_A(r)}{r^q}
 =
 \frac{\arctan z}{z}
 \left(\frac{\sin r}{r}\right)^q
 \frac{1}{\cos r}.
\]
Since $\sin r/r$, $1/\cos r$, and $\arctan z/z$
are smooth even functions near zero, the quotient
$\tau_A(r)/r^q$ can be written as a smooth positive
function of $r^2$ and $r^{2q}$, with value one at the origin.
In particular,
\begin{equation}\label{eq:Atkinson-pole}
 \tau_A(r)
 =r^q\bigl(1+O(r^2)+O(r^{2q})\bigr).
\end{equation}
Differentiating the above representation and using
\eqref{eq:Atkinson-derivative}, we obtain the bounds
in \eqref{eq:profile-bounds} through order $M+1$.
By the symmetry $\pi-\tau_A(\pi-r)=\tau_A(r)$,
the same bounds hold at the opposite pole $-x$.
For convergence analyses of product quadrature rules
based on this transformation, see
\cite{AtkinsonSommariva2005,Sidi2005Atkinson}.

\textbf{Sidi's transformation.} For a real parameter $m>0$, Sidi's $\sin^m$ transformation in 
\cite{Sidi2005Sphere,Sidi2006Surface} is defined by
\begin{equation}\label{eq:Sidi-psi}
 \tau_S(\theta)=\frac{\pi}{c_m}
 \int_0^{\theta/\pi}\sin^m(\pi u)\dd u,
\end{equation}
where $c_m=\int_0^1\sin^m(\pi u)\dd u$.
Differentiating \eqref{eq:Sidi-psi} gives $\tau_S'(\theta)={\sin^m\theta}/{c_m}$, and hence the surface Jacobian is
\begin{equation*}
 J_S(\theta)
 =\frac{\sin^{m-1}\theta\sin\tau_S(\theta)}{c_m}.
\end{equation*}
To examine the behaviour near $r=0$, we let $q=m+1$. A change of variables in \eqref{eq:Sidi-psi} gives
\begin{equation}\label{eq:Sidi-pole}
 \tau_S(r)=\frac{r^q}{c_{q-1}}
 \int_0^1 u^{q-1}
 \left(\frac{\sin(ru)}{ru}\right)^{q-1}\dd u=\frac{r^q}{q c_{q-1}}\bigl(1+O(r^2)\bigr).
\end{equation}
Since $\sin(ru)/(ru)$ is a smooth even function of $r$,
the quotient $\tau_S(r)/r^q$ is a smooth positive function
of $r^2$ near zero. Differentiating the integral
representation in \eqref{eq:Sidi-pole}, we obtain the
bounds in \eqref{eq:profile-bounds} through order $M+1$.
By the symmetry $\pi-\tau_S(\pi-r)=\tau_S(r)$,
the same bounds hold at the opposite pole $-x$.
Since the grading exponent is $q=m+1$, the condition
\eqref{eq:threshold} for Sidi's transformation reads
\begin{equation}\label{eq:Sidi-condition}
 m>\frac{s+1}{2-\alpha}-1.
\end{equation}

Beyond the transformations of Atkinson and Sidi, our analysis also applies to rational sigmoidal transformations \cite{Elliott1998Sigmoidal} and symmetric Korobov transformations \cite{Sidi2009Endpoint}, after rescaling them from $[0,1]$ to $[0,\pi]$. It is straightforward to verify that these
transformations satisfy Assumption~\ref{ass:grading} with the corresponding endpoint grading exponents.

\subsection{Localized cap-area transformation}\label{sec:cap-area}
We propose a novel transformation based on the area of a spherical cap.
It concentrates the cubature points near the singularity and leaves the
sphere outside the cap unchanged. For a prescribed cap radius
$R\in(0,\pi)$, let
\[
 u=\frac{1-\cos\theta}{1-\cos R},
 \qquad 0\le\theta\le R.
\]
The variable $u$ represents the ratio of the area of a cap of radius $\theta$ to that of a cap of radius $R$, both centred at the target $x$. Working in this coordinate allows us to design a localized grading with a particularly simple surface Jacobian. For $q>1$, choose a nonincreasing function
$\chi\in C^\infty([0,1])$ taking values in $[0,1]$, with $\chi=1$ near zero and $\chi=0$ near one.
We define
\[
 \psi_q(u)=u^{1+(q-1)\chi(u)},
 \qquad 0<u\le1,
\]
and let $\psi_q(0)=0$. Logarithmic differentiation gives
\[
 \frac{\psi_q'(u)}{\psi_q(u)}
 =
 \frac{1+(q-1)\chi(u)}{u}
 +(q-1)\chi'(u)\log u.
\]
Since $\chi(u)\ge0$, $\chi'(u)\le0$, and $\log u\le0$, we have $\psi_q'(u)\ge\psi_q(u)/u>0$ for $0<u<1$. The corresponding radial profile is then defined by
\begin{equation}\label{eq:cap-area-profile}
 \tau_C(\theta)=
 \begin{cases}
 \arccos\bigl(1-(1-\cos R)\psi_q(u)\bigr),
 &0\le\theta<R,\\
 \theta,&R\le\theta\le\pi.
 \end{cases}
\end{equation}
As $\psi_q(u)=u$ near $u=1$, the first expression reduces to $\theta$ in a neighbourhood of the cap boundary. Consequently, $\tau_C$ is smooth and strictly increasing on $(0,\pi)$, and the resulting transformation leaves the sphere outside the cap unchanged.

Differentiating \eqref{eq:cap-area-profile} and using \eqref{eq:surface-J}, we obtain the surface Jacobian
\begin{equation}\label{eq:cap-area-J}
 J_C(\theta)
 =\frac{\tau_C'(\theta)\sin\tau_C(\theta)}{\sin\theta}
 =\psi_q'(u),
 \qquad 0<\theta<R,
\end{equation}
with $J_C(\theta)=1$ outside the cap. This simple expression is a consequence of the surface measure taking the form $\dd S=(1-\cos R)\dd u\,\dd\phi$
in the cap-area coordinate.

It remains to verify Assumption~\ref{ass:grading}. Since $\psi_q(u)=u^q$ near zero, we have
\[
 \tau_C(\theta)
 =2\arcsin\left(
 \frac{\sin^q(\theta/2)}{\sin^{q-1}(R/2)}
 \right)
\]
for sufficiently small $\theta$. As in the analysis of
Atkinson's transformation, the quotient
$\tau_C(\theta)/\theta^q$ can be expressed as a smooth
positive function of $\theta^2$ and $\theta^{2q}$ near
the origin. In particular,
\begin{equation*}
 \tau_C(\theta)
 =[2(1-\cos R)]^{(1-q)/2}\theta^q
 \bigl(1+O(\theta^2)+O(\theta^{2q})\bigr).
\end{equation*}
Since $\tau_C(\theta)/\theta^q$ is a smooth function of $\theta^2$ and
$\theta^{2q}$, differentiation gives
the bounds in \eqref{eq:profile-bounds} through order $M+1$.
Since $\tau_C$ is the identity near the opposite pole,
Assumption~\ref{ass:grading} is satisfied, with constants
independent of the target $x$ for fixed $q$, $R$, and $\chi$.

\section{Sharp Sobolev regularity under grading}\label{sec:regularity}
For a profile $\tau$ satisfying Assumption~\ref{ass:grading}, write
$\mathcal{T}_x=\mathcal{T}_{x,\tau}$ and $J_x=J_{x,\tau}$, and let
$G_xf$ be the transformed integrand \eqref{eq:intro-g} for this profile.
At the target $x$, set $G_xf(x)=0$. If $\tau_-$ satisfies
\eqref{eq:profile-bounds} near the antipode $-x$, set $G_xf(-x)=0$;
if the transformation is the identity there, set
$G_xf(-x)=K(x,-x)f(-x)$.

\begin{theorem}[Sobolev regularity]\label{thm:regularity}
Let $s>1$, $0<\alpha<2$, and $M=\lfloor s\rfloor+1$. Suppose Assumption~\ref{ass:kernel} holds and $\tau$ satisfies Assumption~\ref{ass:grading} with exponent $q$ such that
\begin{equation}\label{eq:threshold}
 q(2-\alpha)>s+1.
\end{equation}
Then, for every $f\in\Hs$,
\begin{equation}\label{eq:mapping}
 \sup_{x\in\Sph}\|G_xf\|_{\Hs}\le C\|K\|_{\alpha,M}\|f\|_{\Hs},
\end{equation}
and
\begin{equation}\label{eq:transformed-integral}
 \mathcal Kf(x)=\int_{\Sph}G_xf(z)\dd S(z).
\end{equation}
The constant depends on $s,\alpha$, the fixed chart radii, and the grading profile $\tau$, but not on $x$, $K$, or $f$.
\end{theorem}

\begin{remark}
Theorem~\ref{thm:regularity} therefore applies to the three transformations in Section~\ref{sec:grading}, under \eqref{eq:threshold} for Atkinson's and localized cap-area transformations, and \eqref{eq:Sidi-condition} for Sidi's transformation.
\end{remark}

\subsection{Proof of the Sobolev regularity theorem}
For the local estimates below, we use the equivalence between the spherical harmonic norm \eqref{eq:Hs} and the Sobolev norms of localized chart representations, as recalled in Section~\ref{sec:sobolev}. For $0<\theta<1$ and an open set $U\subset\R^2$, let
\[
 [v]_{\theta,U}^2 =\int_U\int_U\frac{|v(u)-v(w)|^2}{|u-w|^{2+2\theta}}\dd u\dd w.
\]
If $s=k+\theta$, then on a fixed smooth bounded domain $U$,
\begin{equation}\label{eq:local-sobolev-norm}
 \|v\|_{H^s(U)}^2
 \asymp \sum_{|\nu|\le k}\|D^\nu v\|_{L^2(U)}^2 +\sum_{|\nu|=k}[D^\nu v]_{\theta,U}^2;
\end{equation}
see \cite{DiNezzaEtAl2012}. For integer $s=k$, the fractional sum is omitted. 
\begin{lemma}\label{lem:local-composition}
Let $A,B$ be fixed smooth bounded domains. Suppose $F:A\to B$ is a
smooth one-to-one map onto its image, is bi-Lipschitz with uniformly
bounded Jacobians and inverse Jacobians, and has bounded derivatives
through order $k+1$. If $b\in C^{k+1}(\overline A)$ has support in a
fixed compact subset of $A$, then extension by zero gives
\begin{equation}\label{eq:local-composition}
 \|b(v\circ F)\|_{H^s(\R^2)}
 \le C\|b\|_{C^{k+1}(\overline A)}\|v\|_{H^s(B)},
 \qquad s=k+\theta,\quad 0\le\theta<1.
\end{equation}
The constant is uniform over maps with the stated bounds.
\end{lemma}
\begin{proof}
We first prove the estimate for smooth $v$. For integers $k$, the product rule, chain rule, and change of variables give 
\begin{equation}\label{eq:local-integer-estimate}
 \|b(v\circ F)\|_{H^k(A)}\le C\|b\|_{C^k(\overline A)}\|v\|_{H^k(B)}.
\end{equation}
For $0<\theta<1$, we also need to control the fractional seminorms of
$D^\nu[b(v\circ F)]$ for $|\nu|=k$. By the product and chain rules,
each such derivative is a finite sum of terms
$c((D^\mu v)\circ F)$ with $|\mu|\le k$, where $c$ consists of products
of derivatives of $b$ and $F$ of order at most $k$.
The assumed derivative bounds on $F$ through order $k+1$ therefore give
\[
 \|c\|_\infty+\|c\|_{\mathrm{Lip}} \le C\|b\|_{C^{k+1}(\overline A)}.
\]
The identity $c(x)w(x)-c(y)w(y) =c(x)(w(x)-w(y))+(c(x)-c(y))w(y)$
gives
\begin{equation}\label{eq:fractional-product}
 \begin{aligned}
 [cw]_{\theta,A}^2
 &\le 2\|c\|_\infty^2[w]_{\theta,A}^2
       +2\|c\|_{\mathrm{Lip}}^2
       \int_A|w(y)|^2\left(\int_A|x-y|^{-2\theta}\dd x\right)\dd y\\
 &\le 2\|c\|_\infty^2[w]_{\theta,A}^2+2\|c\|_{\mathrm{Lip}}^2\|w\|_{L^2(A)}^2
       \int_{|h|\le\operatorname{diam}A}|h|^{-2\theta}\dd h\\
 &\le C\bigl(\|c\|_\infty^2[w]_{\theta,A}^2
       +\|c\|_{\mathrm{Lip}}^2\|w\|_{L^2(A)}^2\bigr).
 \end{aligned}
\end{equation}
Bi-Lipschitz change of variables and $F(A)\subset B$ give
\begin{equation}\label{eq:fractional-composition}
 [w\circ F]_{\theta,A}^2
 \le C\int_{F(A)}\int_{F(A)}
       \frac{|w(u)-w(v)|^2}{|u-v|^{2+2\theta}}\dd u\dd v
 \le C[w]_{\theta,B}^2,
\end{equation}
and
\begin{equation}\label{eq:L2-composition}
 \|w\circ F\|_{L^2(A)}\le C\|w\|_{L^2(B)}.
\end{equation}
Applying \eqref{eq:fractional-product} to $w=(D^\mu v)\circ F$, then
using \eqref{eq:fractional-composition}, \eqref{eq:L2-composition},
and the bound on $c$, gives
\[
 [c((D^\mu v)\circ F)]_{\theta,A}^2
 \le C\|b\|_{C^{k+1}(\overline A)}^2
       \bigl([D^\mu v]_{\theta,B}^2+\|D^\mu v\|_{L^2(B)}^2\bigr)\le C\|b\|_{C^{k+1}(\overline A)}^2
       \|D^\mu v\|_{H^\theta(B)}^2.
\]
For $|\mu|=k$, \eqref{eq:local-sobolev-norm} bounds the last norm by
$C\|v\|_{H^s(B)}$. For $|\mu|<k$, a bounded $H^1$ extension from $B$
to $\R^2$ and $(1+|\xi|^2)^\theta\le1+|\xi|^2$ yield
\[
 \|D^\mu v\|_{H^\theta(B)}
 \le C\|D^\mu v\|_{H^1(B)}
 \le C\|v\|_{H^s(B)}.
\]
Summing the finitely many derivative terms therefore gives
\begin{equation}\label{eq:local-fractional-estimate}
 \sum_{|\nu|=k}[D^\nu\{b(v\circ F)\}]_{\theta,A}^2
 \le C\|b\|_{C^{k+1}(\overline A)}^2\|v\|_{H^s(B)}^2.
\end{equation}
The fixed support set has distance at least $d>0$ from $\partial A$.
For $g=D^\nu[b(v\circ F)]$, $|\nu|=k$, extended by zero,
\[
 \begin{aligned}
 [g]_{\theta,\R^2}^2
 &=[g]_{\theta,A}^2+2\int_A\int_{\R^2\setminus A}
       \frac{|g(x)-g(y)|^2}{|x-y|^{2+2\theta}}\dd y\dd x\\
 &\le[g]_{\theta,A}^2+2\|g\|_{L^2(A)}^2
       \int_{|h|\ge d}|h|^{-2-2\theta}\dd h\\
 &=[g]_{\theta,A}^2+\frac{2\pi}{\theta d^{2\theta}}\|g\|_{L^2(A)}^2.
 \end{aligned}
\]
Combining this zero-extension bound with \eqref{eq:local-integer-estimate} and \eqref{eq:local-fractional-estimate} proves \eqref{eq:local-composition} by \eqref{eq:local-sobolev-norm}. The result for general $v\in H^s(B)$ follows by density,
using the $L^2$ continuity of the composition operator
established in \eqref{eq:L2-composition}.
\end{proof}

\begin{lemma}\label{lem:weighted-pullback}
Let $F(\zeta)=\tau(|\zeta|)\zeta/|\zeta|$ with \eqref{eq:profile-bounds}, and let $r_0$ be the radius in Assumption~\ref{ass:grading}. Choose $0<R_0<\min\{1,r_0/4\}$ sufficiently small. Suppose $W$ is defined for $0<|\zeta|<r_0$ and satisfies
\begin{equation}\label{eq:scaled-weight}
 \sup_{0<R\le R_0}
 \|R^{-\beta}W(R\,\cdot)\|_{C^{k+1}(\overline A)}\le C_W,
 \qquad A=\{1/4<|v|<4\}.
\end{equation}
Let $\rho_0(\zeta)=\rho(\zeta/R_0)$, where $\rho\in C_c^\infty(\R^2)$ is fixed, radial, one on the unit ball and zero outside the ball of radius two. If
$\delta=\beta+1-s>0$, then
\begin{equation}\label{eq:weighted-pullback}
 \|\rho_0 W(u\circ F)\|_{H^s(\R^2)}
 \le C C_W\frac{R_0^\delta}{\sqrt{1-2^{-2\delta}}}
          \|u\|_{H^s(\R^2)}.
\end{equation}
The product $\rho_0W(u\circ F)$ extends continuously to zero with value zero.
\end{lemma}
\begin{proof}
Let $\eta(v)=\rho(v)-\rho(2v)$ and $R_j=R_0 2^{-j}$, reducing $R_0$
if necessary. For $\zeta\ne0$, we have $\rho_0(\zeta)=\sum_{j=0}^\infty\eta(\zeta/R_j)$.
The rescaled map
\begin{equation}\label{eq:rescaled-map}
 F_R(v)=R^{-q}F(Rv)
\end{equation}
has uniformly bounded derivatives through order $k+1$ by
\eqref{eq:profile-bounds}. The identity
\[
 R^{2q}|F_R(r\omega)-F_R(t\omega')|^2
 =|\tau(Rr)-\tau(Rt)|^2+\tau(Rr)\tau(Rt)|\omega-\omega'|^2
\]
and \eqref{eq:profile-bounds} give uniform
bi-Lipschitz bounds on $A$. All images $F_R(A)$ lie in a fixed annulus $B$.
For each dyadic term, let 
\[h_j(\zeta)=\eta(\zeta/R_j)W(\zeta)u(F(\zeta)).\]
Its rescaling is $R_j^{-\beta}h_j(R_jv)
 =\eta(v)R_j^{-\beta}W(R_jv)u(R_j^q F_{R_j}(v))$.
We then have
\begin{equation}\label{eq:Vj}
 R_j^{-\beta}\|h_j(R_j\,\cdot)\|_{H^s(\R^2)}
 \le CC_W\|u(R_j^q\,\cdot)\|_{H^s(B)}
 \le CC_W\|u\|_{H^s(\R^2)}.
\end{equation}
Indeed, applying Lemma~\ref{lem:local-composition} with $b=\eta R_j^{-\beta}W(R_j\,\cdot)$, $F=F_{R_j}$, and $v=u(R_j^q\,\cdot)$, and using \eqref{eq:scaled-weight}, gives the first inequality. For the second inequality, let $t=R_j^q\le1$ and $s=k+\theta$, $0\le\theta<1$. We bound each term in the local Sobolev norm \eqref{eq:local-sobolev-norm}. The zeroth-order term satisfies $\|u(t\,\cdot)\|_{L^2(B)}\le |B|^{1/2}\|u\|_\infty$.
For $1\le|\nu|\le k$, the chain rule and change of variables give
\[
 \|D^\nu(u(t\,\cdot))\|_{L^2(B)}
 =t^{|\nu|-1}\|D^\nu u\|_{L^2(tB)}
 \le\|D^\nu u\|_{L^2(\R^2)}.
\]
For $0<\theta<1$ and $|\nu|=k$, the same change of variables in the double integral gives
\[
 [D^\nu(u(t\,\cdot))]_{\theta,B}
 =t^{s-1}[D^\nu u]_{\theta,tB}
 \le[D^\nu u]_{\theta,\R^2}.
\]
Here $t^{|\nu|-1}\le1$ and $t^{s-1}\le1$, and each integral over $tB$ is bounded by the corresponding integral over $\R^2$. Combining these estimates in \eqref{eq:local-sobolev-norm} and using the Sobolev embedding $H^s(\R^2)\hookrightarrow L^\infty(\R^2)$, we obtain
\[
 \|u(t\,\cdot)\|_{H^s(B)}^2\le C\left(
       \|u\|_\infty^2
       +\sum_{1\le|\nu|\le k}\|D^\nu u\|_{L^2(\R^2)}^2
       +\sum_{|\nu|=k}[D^\nu u]_{\theta,\R^2}^2
       \right)\le C\|u\|_{H^s(\R^2)}^2,
\]
with the fractional sum omitted when $\theta=0$. The constant is independent of $t$, and hence of $j$, proving the second inequality in \eqref{eq:Vj}.

For $0<R\le1$,  Fourier scaling gives
\begin{equation}\label{eq:fourier-scaling}
 \left\|u\left({\cdot}/{R}\right)\right\|_{H^s(\R^2)}^2
 =R^2\int_{\R^2}(1+|\xi|^2/R^2)^s|\widehat u(\xi)|^2\dd\xi
\le R^{2-2s}\|u\|_{H^s(\R^2)}^2.
\end{equation}
Combining \eqref{eq:Vj} and \eqref{eq:fourier-scaling} gives
\(\|h_j\|_{H^s}\le CC_WR_j^\delta\|u\|_{H^s}.\)
We use finite overlap before summing these bounds. The supports of $h_j$ and of all their derivatives lie in $\{R_j/2\le|\zeta|\le2R_j\}$. At each point only a uniformly bounded number of these annuli meet. For any finite index set $I$, bounded overlap and Cauchy--Schwarz give 
\begin{equation}\label{eq:finite-overlap}
 \|\sum_{j\in I}h_j\|_{H^s(\R^2)}^2
 \le C_s\sum_{j\in I}\|h_j\|_{H^s(\R^2)}^2
\end{equation}
for integer $s$. If $s=k+\theta$ with $0<\theta<1$, we take $|\nu|=k$. Only a uniformly bounded number of indices contribute to
the difference at a pair $(\zeta,v)$, by finite overlap at both points. Thus
\[
 \left|\sum_{j\in I}[D^\nu h_j(\zeta)-D^\nu h_j(v)]\right|^2
 \le C\sum_{j\in I}|D^\nu h_j(\zeta)-D^\nu h_j(v)|^2.
\]
Integrating this inequality against
$|\zeta-v|^{-2-2\theta}\dd\zeta\dd v$ and combining it with the
integer derivative bounds proves the same estimate \eqref{eq:finite-overlap} for fractional $s$.

Since $\delta>0$,
\begin{equation*}
 \sum_{j=0}^\infty R_j^{2\delta}=\frac{R_0^{2\delta}}{1-2^{-2\delta}}<\infty.
\end{equation*}
Estimates \eqref{eq:Vj}, \eqref{eq:fourier-scaling}, and
\eqref{eq:finite-overlap} therefore imply that $\sum_{j=0}^\infty h_j$ converges
in $H^s(\R^2)$. For each $\zeta\ne0$, only finitely many terms are
nonzero, and their sum is $\rho_0W(u\circ F)$. Moreover,
\eqref{eq:scaled-weight} gives $W(\zeta)=O(|\zeta|^\beta)$ near zero.
Since $\beta>s-1>0$ and $u$ is bounded, this product tends to zero as
$\zeta\to0$. Its extension with value zero at the origin is therefore
the continuous representative of the $H^s$ limit. Finally,
\eqref{eq:finite-overlap} and the geometric sum above give
\eqref{eq:weighted-pullback}.
\end{proof}

Together with Lemmata \ref{lem:local-composition} and \ref{lem:weighted-pullback}, we are ready to prove Theorem~\ref{thm:regularity}.

\begin{proof}[Proof of Theorem~\ref{thm:regularity}]
We first take $f$ smooth. With the cutoff $\rho$ from Lemma~\ref{lem:weighted-pullback}, we define
\[
 u_x(\zeta):=\rho(\zeta/r_0)f(\exp_x\zeta),\qquad \zeta\in\R^2
\]
in normal coordinates at $x$. Then 
\begin{equation*}
u_x(\zeta)=\begin{cases}
f(\exp_x\zeta), &|\zeta|\le r_0,\\
0, & |\zeta|\ge2r_0.
\end{cases}
\end{equation*}
The Sobolev norm equivalence gives $\|u_x\|_{H^s(\R^2)}\le C_s\|f\|_{\Hs}$ uniformly in $x$ and the tangent frame. Locally,
\[
 G_xf(\exp_x\zeta)=W_x(\zeta)u_x(F(\zeta)),
\]
where $W_x(\zeta)=K(x,\exp_x F(\zeta))J_x(\exp_x\zeta)$. We verify \eqref{eq:scaled-weight} with $\beta=q(2-\alpha)-2$.
For $r=|v|$, formula \eqref{eq:surface-J} gives
\begin{equation}\label{eq:scaled-J}
 R^{-(2q-2)}J_x(\exp_x(Rv))
 =\big[R^{1-q}\tau'(Rr)\big]
   \frac{\tau(Rr)}{R^q r}
   \frac{\sin\tau(Rr)/\tau(Rr)}{\sin(Rr)/(Rr)}.
\end{equation}
The right side has a uniform $C^{k+1}(A)$ bound by
\eqref{eq:profile-bounds} through order $k+2=M+1$.
With $F_R$ from \eqref{eq:rescaled-map},
\begin{equation*}
 R^{-\beta}W_x(Rv)
 =\big[R^{q\alpha}K(x,\exp_x(R^qF_R(v)))\big]
   R^{-(2q-2)}J_x(\exp_x(Rv)).
\end{equation*}
The kernel derivative bounds \eqref{eq:kernel-derivatives}, applied at
scale $R^q$, together with \eqref{eq:scaled-J} and the uniform bounds
on $F_R$ proved in Lemma~\ref{lem:weighted-pullback}, give
\[
 \sup_{x,\,R}\|R^{-\beta}W_x(R\,\cdot)\|_{C^{k+1}(A)}
 \le C\|K\|_{\alpha,M}.
\]
Condition \eqref{eq:threshold} says exactly $\beta>s-1$.
Lemma~\ref{lem:weighted-pullback} therefore bounds the contribution near
the target $x$ by the right-hand side of \eqref{eq:mapping}.

Near $-x$, the kernel is $C^M$. If $\tau_-$ satisfies
\eqref{eq:profile-bounds}, the same argument applies with $\beta=2q-2$, since $2q-s-1=\big[q(2-\alpha)-s-1\big]+q\alpha>0$.
If the transformation is the identity near $-x$,
Lemma~\ref{lem:local-composition} applies there directly. Away from both
poles, $\mathcal T_x$ is a smooth diffeomorphism with uniformly bounded
derivatives and inverse derivatives, and its image stays away from $x$.
Lemma~\ref{lem:local-composition} therefore gives the remaining bounds
in finitely many charts. Combining these local bounds using the Sobolev
norm equivalence in Section~\ref{sec:sobolev} proves \eqref{eq:mapping}
for smooth $f$.

For $f\in\Hs$, we take smooth $f_n\to f$ in $\Hs$. Estimate \eqref{eq:mapping} makes $G_xf_n$ Cauchy in $\Hs$, uniformly
in $x$. Since the transformed weight is bounded, \eqref{eq:embedding}
also gives uniform convergence to $G_xf$. Passing to the limit proves
\eqref{eq:mapping} for $f$.

Finally, recall that $K(x,\cdot)f$ is absolutely integrable by \eqref{eq:kernel-pointwise}, and $G_xf\in\Hs\subset C(\Sph)$.
Changing variables on $\Sph\setminus\{x,-x\}$, where $\mathcal T_x$
is a smooth diffeomorphism, proves \eqref{eq:transformed-integral}.
\end{proof}

\subsection{Sharpness of the grading condition}\label{sec:sharpness}
We give a counterexample showing that the strict inequality in
\eqref{eq:threshold} is necessary for the kernel class in
Assumption~\ref{ass:kernel}.

\begin{theorem}[Sharp grading threshold]\label{thm:sharp-threshold}
Let $s>1$, $0<\alpha<2$, and let $\tau$ satisfy
Assumption~\ref{ass:grading} with exponent $q$. If $q(2-\alpha)\le s+1$,
there is a kernel $K$ satisfying Assumption~\ref{ass:kernel} such that
\begin{equation}\label{eq:sharp-counterexample}
 G_{x_0}f\notin H^s(\Sph),
\end{equation}
where $ f\equiv1$ and $ x_0=(0,0,1)$. 
The kernel can be chosen smooth off the diagonal, strictly positive,
and uniformly comparable to $|p-y|^{-\alpha}$ for $p\ne y$.
\end{theorem}
\begin{proof}
Choose $\chi\in C^\infty([0,\infty))$, with $0\le\chi\le1$, such that
\[
 \chi(\rho)=
 \begin{cases}
 1,&0\le\rho\le r_0/4,\\
 0,&\rho\ge r_0/2.
 \end{cases}
\]
For $p\ne y$, we define
\[
 K(p,y)=
 \begin{cases}
 \chi(\rho)\rho^{-\alpha}\bigl(1+\tfrac12\omega_1^{k+1}\bigr)
       +1-\chi(\rho),&0<\rho<r_0,\\
 1,&\rho\ge r_0,
 \end{cases}
 \qquad \rho=d(p,y).
\]
In the first case, $y=\exp_p(\rho\omega)$ and $\omega_1$ is the first
Cartesian component of the unit tangent vector $\omega$. The angular factor lies in
$[1/2,3/2]$, so $K$ is positive and comparable to $|p-y|^{-\alpha}$.
It is smooth off the diagonal and satisfies Assumption~\ref{ass:kernel}
by \eqref{eq:directional}. At $x_0$ and for small $r$, the transformed integrand is
\begin{equation*}
 g(r,\phi):=G_{x_0}1(\exp_{x_0}(r\omega))
 =w(r)\bigl(1+\tfrac12\cos^{k+1}\phi\bigr),
\end{equation*}
where $ w(r)=\tau(r)^{-\alpha}\tau'(r){\sin\tau(r)}/{\sin r}$.
With $\beta=q(2-\alpha)-2$, \eqref{eq:profile-bounds} and
\eqref{eq:scaled-J} show that $R^{-\beta}w(R\,\cdot)$ is uniformly
bounded in $C^{k+1}([1,2])$ and bounded below by a positive constant.
On $A=\{v\in\mathbb{R}^2:1<|v|<2\}$, the highest-order terms in
\eqref{eq:local-sobolev-norm} for $R^{-\beta}g(R\,\cdot)$ therefore
have a uniform positive lower bound. Indeed, otherwise compactness in
$C^k(\overline A)$ and Fatou's lemma would give a limit with vanishing
highest-order terms, hence a polynomial of degree at most $k$.
This is impossible because its angular Fourier mode of order $k+1$
is nonzero.

Let $s=k+\theta$ for $0\le\theta<1$. Scaling back to the disjoint annuli
$A_j=R_jA$, where $R_j=R_0 2^{-j}$, gives
\begin{equation*}
  cR_j^{2(\beta+1-s)}\le \begin{cases}
 \displaystyle\sum_{|\nu|=k}\|D^\nu g\|_{L^2(A_j)}^2,
       &\theta=0,\\[3pt]
 \displaystyle\sum_{|\nu|=k}[D^\nu g]_{\theta,A_j}^2,
       &0<\theta<1.
 \end{cases}
\end{equation*}
Since $\beta+1-s\le0$, summing these lower bounds over $j$ diverges,
whereas \eqref{eq:local-sobolev-norm} bounds their sum whenever
$g\in H^s$ near zero. Thus \eqref{eq:sharp-counterexample} follows.
Together with Theorem~\ref{thm:regularity}, this proves that
\eqref{eq:threshold} is necessary and sufficient for \eqref{eq:mapping}
to hold for every kernel in Assumption~\ref{ass:kernel}.
\end{proof}

\section{Graded cubature with QMC designs}\label{sec:cubature}
We apply the graded rule \eqref{eq:graded-rule} with a profile satisfying
Assumption~\ref{ass:grading} and the notation of
Section~\ref{sec:regularity}. For the equal-weight rule \eqref{equ:quad},
its evaluation takes the form
\[
 \mathcal K_N^{\rm gr}f(x)
 =\frac{4\pi}{N}\sum_{i=1}^N
 K(x,\mathcal T_xz_{N,i})f(\mathcal T_xz_{N,i})J_x(z_{N,i}).
\]
At a node equal to $x$ or $-x$, the integrand value in this sum is
\[
 G_xf(x)=0,\qquad
 G_xf(-x)=
 \begin{cases}
 0,&\tau_-\text{ satisfies }\eqref{eq:profile-bounds},\\
 K(x,-x)f(-x),&\mathcal T_x=\mathrm{id}\text{ near }-x.
 \end{cases}
\]
The reference nodes $z_{N,i}$ are shared
by all targets, while the density is evaluated at the transformed nodes
$\mathcal T_xz_{N,i}$ and the Jacobian supplies the change-of-variables
weights. For the localized cap-area transformation
\eqref{eq:cap-area-profile}, nodes outside the cap $d(x,z)<R$ stay fixed
and have $J_x=1$ by \eqref{eq:cap-area-J}.

Theorem~\ref{thm:regularity} places the transformed integrand in the
Sobolev space used to measure the cubature error. The resulting estimate
is uniform in the target.

\begin{theorem}[Cubature error]\label{thm:cubature}
Under the hypotheses of Theorem~\ref{thm:regularity}, let $Q_N$ be an
$H^s$ QMC-design family satisfying \eqref{eq:intro-qmc}. Then
\begin{equation}\label{eq:QMC-error}
 \|\mathcal Kf-\mathcal K_N^{\rm gr}f\|_\infty
 \le CC_Q\|K\|_{\alpha,M}N^{-s/2}\|f\|_{\Hs}.
\end{equation}
\end{theorem}
\begin{proof}
The integral identity \eqref{eq:transformed-integral} and \eqref{eq:intro-qmc} give
\[
 |\mathcal Kf(x)-Q_N(G_xf)|\le C_QN^{-s/2}\|G_xf\|_{\Hs}.
\]
Taking the supremum over $x$ and applying \eqref{eq:mapping} gives
\eqref{eq:QMC-error}.
\end{proof}

The error estimate also applies at lower Sobolev orders. Let
$f\in H^s(\Sph)$ and suppose
$\text{wce}(Q_N;H^r)\le C_rN^{-r/2}$.
If Assumptions~\ref{ass:kernel} and \ref{ass:grading} hold at order $r$,
with $M=\lfloor r\rfloor+1$, $1<r\le s$, and $r<q(2-\alpha)-1$, then Theorem~\ref{thm:cubature} at order $r$ and
$\|f\|_{H^r}\le\|f\|_{H^s}$ give
\begin{equation}\label{eq:independent-error}
 \|\mathcal Kf-\mathcal K_N^{\rm gr}f\|_\infty
 \le C C_r\|K\|_{\alpha,\lfloor r\rfloor+1}N^{-r/2}\|f\|_{H^s(\Sph)}.
\end{equation}
Thus the order $r$ is limited by the density regularity, the grading
exponent, and the Sobolev worst-case error of the point family.
For an equal-weight family, its QMC strength is
\begin{equation*}
 \sigma=\sup\{t>1:\ \sup_N N^{t/2}\text{wce}(Q_N;H^t)<\infty\},
\end{equation*}
assuming this set is nonempty. By \cite[Theorem 4]{MR3246811}, the
QMC-design bound \eqref{eq:intro-qmc} holds at every order
$1<r<\sigma$. Consequently, under the kernel and profile assumptions
at order $r$, \eqref{eq:independent-error} holds whenever
\[
 1<r<\min\{\sigma,q(2-\alpha)-1\},\qquad r\le s.
\]
The endpoint $r=\sigma\le s$ is also available if the bound in
\eqref{eq:intro-qmc} holds there; the grading inequality remains strict.

For a prescribed order $r$ supported by the density and the point
family, it therefore suffices to choose
$q>(r+1)/(2-\alpha)$. For example, when $\alpha=1$ and $f\in H^2(\Sph)$,
an $H^2$ QMC-design family gives an error of order $N^{-1}$ for every
$q>3$. For Sidi's transformation, this condition is $m>2$ by
\eqref{eq:Sidi-condition}. 

\section{Numerical experiments}\label{sec:numerics}
All rules have weights $4\pi/N$, with $N=(t+1)^2$ and $t=10,14,20,28,40,50,64,80,100$. We use spherical $t$-designs (STD) \cite{delsarte1991geometriae}, maximal-determinant (MD) \cite{SloanWomersley2004Extremal}, minimal-energy (ME) \cite{RakhmanovSaffZhou1994}, and equal-area (EA) points \cite{Leopardi2006EqualArea}. We compare Atkinson, Sidi, and localized cap-area transformations, with the same angular exponent $q$; the Sidi parameter is $m=q-1$. For \eqref{eq:cap-area-profile}, we fix $R=2.888$ and use the same smooth
nonincreasing cutoff throughout, equal to one for $u\le0.022212$ and zero
for $u\ge0.97168$.

\subsection{Graded cubature for a direction-dependent kernel}\label{sec:num-cubature}
To examine the effect of grading on a direction-dependent singularity,
we fix the target $x=e_3=(0,0,1)$ and consider the kernel
\begin{equation}\label{eq:numerical-model}
 K(e_3,\exp_{e_3}(r\omega))
 =\frac{\chi(r)}{r}\bigl(1+0.3\cos(4\phi)\bigr)+1-\chi(r),
\end{equation}
where $\omega=(\cos\phi,\sin\phi,0)$ and $\chi$ is a smooth cutoff
function taking values in $[0,1]$, with $\chi(r)=1$ for $r\le0.8$ and
$\chi(r)=0$ for $r\ge1.3$. The kernel has a weak singularity of order
$\alpha=1$ with a direction-dependent leading coefficient.

We first take the constant density $f\equiv1$ to isolate the effect of
grading from the regularity of the density. For $s=2$, the critical
grading exponent in \eqref{eq:threshold} is $q_c=3$. We compare the
three transformations with $q=2.5$, $3$, and $4$, respectively below,
at, and above this threshold.

\begin{figure}[htbp]
 \centering
 \includegraphics[width=\linewidth]{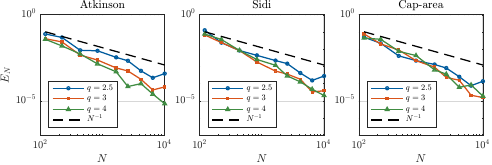}
 \caption{STD cubature errors for the direction-dependent kernel
 \eqref{eq:numerical-model} with $q=2.5,3,4$, under Atkinson, Sidi,
 and localized cap-area transformations.}
 \label{fig:regularization}
\end{figure}

Figure~\ref{fig:regularization} compares the resulting STD cubature errors for the three values of $q$ under each transformation. Let $E_N$ denote the largest absolute cubature error over eight fixed rotations of each point set, which reduces the influence of a particular orientation relative to the singularity. The reference slope $N^{-1}$ in Figure~\ref{fig:regularization} is the rate guaranteed by Theorem~\ref{thm:cubature} for $s=2$ and $q=4$.
Reference values are obtained by analytically integrating out the angular dependence and numerically evaluating the remaining radial
integral. In particular, for $f\equiv1$,
\[
 (\mathcal K1)(e_3)
 =2\pi\int_0^\pi
 \left(\frac{\chi(r)}{r}+1-\chi(r)\right)\sin r\,\dd r,
\]
since the angular mode $\cos(4\phi)$ integrates to zero.

\begin{figure}[htbp]
 \centering
 \includegraphics[width=\linewidth]{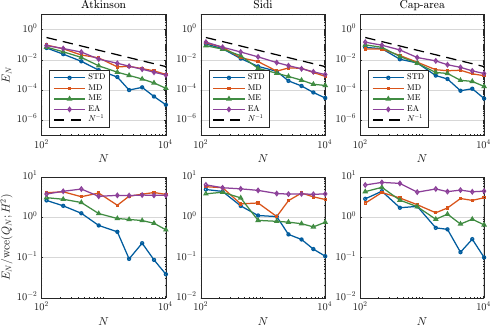}
 \caption{Finite-regularity density under Atkinson, Sidi, and localized cap-area
 transformations.}
 \label{fig:cubature-transfer}
\end{figure}

For the same kernel \eqref{eq:numerical-model}, we now consider a density of finite Sobolev regularity,
\[
 f(y)=1+\frac14|y+e_3|^{5/4}.
\]
This function belongs to $H^r(\Sph)$ for $r<9/4$ but is not smooth at the antipode $-e_3$. In particular, $f\in H^2(\Sph)$, so we apply all three transformations with grading exponent $q=4$ to examine how the worst-case error bounds of the underlying cubature rules translate into accuracy for weakly singular integrals. For the equal-weight rule \eqref{equ:quad}, Parseval's identity and the addition theorem for spherical harmonics
yield
\begin{equation*}
 \mathrm{wce}(Q_N;H^s)^2
 =
 \frac{4\pi}{N^2}\sum_{\ell=1}^{\infty}
 \frac{2\ell+1}{(1+\ell(\ell+1))^s}
 \sum_{i,j=1}^N P_\ell(z_{N,i}\cdot z_{N,j}),
\end{equation*}
where $P_\ell$ is the Legendre polynomial of degree $\ell$, normalized
by $P_\ell(1)=1$. We sum the kernel series through degree 12000 and
interpolate it on 65537 equally spaced polar angles. The remaining
diagonal contribution
$\frac{4\pi}{N}\sum_{\ell>12000}(2\ell+1)/(1+\ell(\ell+1))^2$
is added to the squared error.
Over all point families and values of $N$ tested, the computed values satisfy
\[
 2.72\le N\,\mathrm{wce}(Q_N;H^2)\le2.97,
\]
which is consistent with the $O(N^{-1})$ worst-case
error rate of QMC designs for $H^2(\Sph)$.
Figure~\ref{fig:cubature-transfer} displays the cubature
errors $E_N$ and the normalized errors
$E_N/\mathrm{wce}(Q_N;H^2)$.

\subsection{Laplace double layers on a nonconvex surface}
\begin{figure}[htbp]
 \centering
 \includegraphics[width=\linewidth]{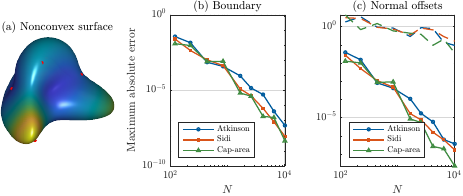}
 \caption{Laplace double-layer evaluation with the three transformations.
 Left: surface and boundary targets. Middle: boundary errors.
 Right: close-evaluation errors with (solid) and without (dashed)
 constant subtraction.}
 \label{fig:surface-double-layer}
\end{figure}

We apply graded QMC cubature to a Laplace double-layer integral pulled back from a smooth nonconvex surface. Recall that $P_j$ denotes the Legendre
polynomial of degree $j$. For $w\in\Sph$, let
\[\rho(w)=\exp\big(0.45P_3(w_3)+0.28(w_1^3-3w_1w_2^2)
                          +0.15P_5(w_2)\big)\]
                          and
                          \[
 \Phi(w)=\bigl(1.4\rho(w)w_1,\,\rho(w)w_2+0.22[1.4\rho(w)w_1]^2,\,0.85\rho(w)w_3
          +0.18[1.4\rho(w)w_1]\rho(w)w_2\bigr).\]
Figure~\ref{fig:surface-double-layer}(a) shows $\Phi(\Sph)$ and the six
boundary targets.  With $n(y)$ the outward unit normal at $\Phi(y)$,
the pulled-back boundary kernel is
\begin{equation}\label{eq:double-layer}
 K_{\rm DL}(x,y)=
 \frac{n(y)\cdot[\Phi(x)-\Phi(y)]}{4\pi|\Phi(x)-\Phi(y)|^3}j_\Phi(y).
\end{equation}
Its numerator is $O(d(x,y)^2)$, giving $\alpha=1$. We prescribe
$\mu=u|_{\Phi(\Sph)}$, where
\[
 u(Y)=Y_1+0.2(Y_2^2-Y_3^2)+0.1Y_1Y_2.
\]
For this smooth density, $q=4$ satisfies \eqref{eq:threshold} with
$s=2$, so Theorem~\ref{thm:cubature} gives an $N^{-1}$ QMC-design
error bound for the boundary integral.

Let $h=\mu\circ\Phi$ and $X_\varepsilon=\Phi(x)+\varepsilon n(x)$.
Using the classical constant-subtraction identity
(see, e.g., \cite{KhatriEtAl2020Close}), we write
\begin{equation}\label{eq:subtracted-DL}
 D\mu(X_\varepsilon)
 =c_\varepsilon h(x)+\int_{\Sph}
 \frac{n(y)\cdot[X_\varepsilon-\Phi(y)]}
 {4\pi|X_\varepsilon-\Phi(y)|^3}
 [h(y)-h(x)]j_\Phi(y)\dd S(y),
\end{equation}
where $c_\varepsilon=-1$ for $\varepsilon<0$, $c_0=-1/2$, and
$c_\varepsilon=0$ for $\varepsilon>0$. At $\varepsilon=0$,
$D\mu(X_0)$ denotes the boundary integral with kernel
\eqref{eq:double-layer}.
Since $u$ is harmonic, Green's identity gives the reference values
\[
 D\mu(X_\varepsilon)
 =S(\partial_nu)(X_\varepsilon)+c_\varepsilon u(X_\varepsilon),
\]
where $S$ is the Laplace single-layer operator. We compute $S(\partial_nu)$
using Atkinson's transformation with $q=5$, followed by Gauss--trapezoidal
product integration at resolution $1536\times3072$; halving the resolution
changes the reference values by less than $10^{-13}$.

We use all three transformations with $q=4$ and
$\varepsilon=0,-10^{-2},-10^{-4},-10^{-6},10^{-4}$.
Figure~\ref{fig:surface-double-layer}(b) shows the boundary errors
at $\varepsilon=0$, maximized over the six boundary targets.
In panel (c), the solid curves use \eqref{eq:subtracted-DL}, with
$c_\varepsilon h(x)$ added exactly; their errors are maximized over
the six targets and all four nonzero offsets. The dashed curves fix
$\varepsilon=-10^{-6}$ and apply the same graded cubature directly to
the double-layer kernel multiplied by $h(y)$, without subtracting
$h(x)$. These errors are maximized over the six targets only.
At $N=10201$ and $q=4$, the STD errors are smallest for Atkinson's
transformation in the model integrals and for the localized cap-area
transformation in the boundary and close evaluations.

\bibliographystyle{siamplain}
{\small \bibliography{myref}}

\begin{thebibliography}{10}

\bibitem{AnChen2016}
{\sc C.~An and S.~Chen}, {\em Numerical integration over the unit sphere by
  using spherical {$t$}-design}, arXiv preprint arXiv:1611.02785,  (2016).

\bibitem{MR5025409}
{\sc C.~An and H.-N. Wu}, {\em Spherical configurations and quadrature methods
  for integral equations of the second kind}, SIAM J. Numer. Anal., 64 (2026),
  pp.~148--169, \url{https://doi.org/10.1137/24M1688370}.

\bibitem{Atkinson1982Sphere}
{\sc K.~Atkinson}, {\em Numerical integration on the sphere}, J. Austral. Math.
  Soc. Ser. B, 23 (1982), pp.~332--347,
  \url{https://doi.org/10.1017/S0334270000000278}.

\bibitem{Atkinson2004Surface}
{\sc K.~Atkinson}, {\em Quadrature of singular integrands over surfaces},
  Electron. Trans. Numer. Anal., 17 (2004), pp.~133--150.

\bibitem{MR2934227}
{\sc K.~Atkinson and W.~Han}, {\em Spherical harmonics and approximations on
  the unit sphere: An introduction}, vol.~2044 of Lecture Notes in Mathematics,
  Springer, Heidelberg, 2012, \url{https://doi.org/10.1007/978-3-642-25983-8}.

\bibitem{AtkinsonSommariva2005}
{\sc K.~Atkinson and A.~Sommariva}, {\em Quadrature over the sphere}, Electron.
  Trans. Numer. Anal., 20 (2005), pp.~104--118.

\bibitem{bondarenko2013optimal}
{\sc A.~Bondarenko, D.~Radchenko, and M.~Viazovska}, {\em Optimal asymptotic
  bounds for spherical designs}, Ann. of Math. (2), 178 (2013), pp.~443--452,
  \url{https://doi.org/10.4007/annals.2013.178.2.2}.

\bibitem{BrauchartHesse2007}
{\sc J.~S. Brauchart and K.~Hesse}, {\em Numerical integration over spheres of
  arbitrary dimension}, Constr. Approx., 25 (2007), pp.~41--71,
  \url{https://doi.org/10.1007/s00365-006-0629-4}.

\bibitem{MR3246811}
{\sc J.~S. Brauchart, E.~B. Saff, I.~H. Sloan, and R.~S. Womersley}, {\em Q{MC}
  designs: optimal order quasi {M}onte {C}arlo integration schemes on the
  sphere}, Math. Comp., 83 (2014), pp.~2821--2851,
  \url{https://doi.org/10.1090/S0025-5718-2014-02839-1}.

\bibitem{DaiXu2013}
{\sc F.~Dai and Y.~Xu}, {\em Approximation Theory and Harmonic Analysis on
  Spheres and Balls}, Springer Monographs in Mathematics, Springer, New York,
  2013, \url{https://doi.org/10.1007/978-1-4614-6660-4}.

\bibitem{delsarte1991geometriae}
{\sc P.~Delsarte, J.-M. Goethals, and J.~J. Seidel}, {\em Spherical codes and
  designs}, Geom. Dedicata, 6 (1977), pp.~363--388,
  \url{https://doi.org/10.1007/bf03187604}.

\bibitem{DiNezzaEtAl2012}
{\sc E.~Di~Nezza, G.~Palatucci, and E.~Valdinoci}, {\em Hitchhiker's guide to
  the fractional {Sobolev} spaces}, Bull. Sci. Math., 136 (2012), pp.~521--573,
  \url{https://doi.org/10.1016/j.bulsci.2011.12.004}.

\bibitem{Duffy1982}
{\sc M.~G. Duffy}, {\em Quadrature over a pyramid or cube of integrands with a
  singularity at a vertex}, SIAM J. Numer. Anal., 19 (1982), pp.~1260--1262,
  \url{https://doi.org/10.1137/0719090}.

\bibitem{Elliott1998Sigmoidal}
{\sc D.~Elliott}, {\em Sigmoidal transformations and the trapezoidal rule}, J.
  Austral. Math. Soc. Ser. B, 40 (1998), pp.~E77--E137,
  \url{https://doi.org/10.21914/anziamj.v40i0.455}.

\bibitem{FariaPerezBonnet2021}
{\sc L.~M. Faria, C.~P{\'e}rez-Arancibia, and M.~Bonnet}, {\em General-purpose
  kernel regularization of boundary integral equations via density
  interpolation}, Comput. Methods Appl. Mech. Engrg., 378 (2021), p.~113703,
  \url{https://doi.org/10.1016/j.cma.2021.113703}.

\bibitem{GrahamSloan2002}
{\sc I.~G. Graham and I.~H. Sloan}, {\em Fully discrete spectral boundary
  integral methods for {Helmholtz} problems on smooth closed surfaces in
  {$\mathbb R^3$}}, Numer. Math., 92 (2002), pp.~289--323,
  \url{https://doi.org/10.1007/s002110100343}.

\bibitem{HesseSloan2005Worst}
{\sc K.~Hesse and I.~H. Sloan}, {\em Worst-case errors in a {Sobolev} space
  setting for cubature over the sphere {$S^2$}}, Bull. Austral. Math. Soc., 71
  (2005), pp.~81--105, \url{https://doi.org/10.1017/S0004972700038041}.

\bibitem{HesseSloan2006Sobolev}
{\sc K.~Hesse and I.~H. Sloan}, {\em Cubature over the sphere {$S^2$} in
  {Sobolev} spaces of arbitrary order}, J. Approx. Theory, 141 (2006),
  pp.~118--133, \url{https://doi.org/10.1016/j.jat.2006.01.004}.

\bibitem{JinZhou2015}
{\sc B.~Jin and Z.~Zhou}, {\em A finite element method with singularity
  reconstruction for fractional boundary value problems}, ESAIM Math. Model.
  Numer. Anal., 49 (2015), pp.~1261--1283,
  \url{https://doi.org/10.1051/m2an/2015010}.

\bibitem{KhatriEtAl2020Close}
{\sc S.~Khatri, A.~D. Kim, R.~Cortez, and C.~Carvalho}, {\em Close evaluation
  of layer potentials in three dimensions}, Journal of Computational Physics,
  423 (2020), p.~109798, \url{https://doi.org/10.1016/j.jcp.2020.109798}.

\bibitem{KlocknerEtAl2013}
{\sc A.~Kl{\"o}ckner, A.~Barnett, L.~Greengard, and M.~O'Neil}, {\em Quadrature
  by expansion: A new method for the evaluation of layer potentials}, J.
  Comput. Phys., 252 (2013), pp.~332--349,
  \url{https://doi.org/10.1016/j.jcp.2013.06.027}.

\bibitem{Leopardi2006EqualArea}
{\sc P.~Leopardi}, {\em A partition of the unit sphere into regions of equal
  area and small diameter}, Electron. Trans. Numer. Anal., 25 (2006),
  pp.~309--327.

\bibitem{RakhmanovSaffZhou1994}
{\sc E.~A. Rakhmanov, E.~B. Saff, and Y.~M. Zhou}, {\em Minimal discrete energy
  on the sphere}, Math. Res. Lett., 1 (1994), pp.~647--662,
  \url{https://doi.org/10.4310/MRL.1994.v1.n6.a3}.

\bibitem{Sidi2005Atkinson}
{\sc A.~Sidi}, {\em Analysis of {Atkinson}'s variable transformation for
  numerical integration over smooth surfaces in {$\mathbb R^3$}}, Numer. Math.,
  100 (2005), pp.~519--536, \url{https://doi.org/10.1007/s00211-005-0605-6}.

\bibitem{Sidi2005Sphere}
{\sc A.~Sidi}, {\em Application of class {$\mathcal S_m$} variable
  transformations to numerical integration over surfaces of spheres}, J.
  Comput. Appl. Math., 184 (2005), pp.~475--492,
  \url{https://doi.org/10.1016/j.cam.2005.01.022}.

\bibitem{Sidi2006Surface}
{\sc A.~Sidi}, {\em Numerical integration over smooth surfaces in {$\mathbb
  R^3$} via class {$\mathcal S_m$} variable transformations. {Part II}:
  Singular integrands}, Appl. Math. Comput., 181 (2006), pp.~291--309,
  \url{https://doi.org/10.1016/j.amc.2006.01.035}.

\bibitem{Sidi2009Endpoint}
{\sc A.~Sidi}, {\em Variable transformations and {Gauss--Legendre} quadrature
  for integrals with endpoint singularities}, Math. Comp., 78 (2009),
  pp.~1593--1612, \url{https://doi.org/10.1090/S0025-5718-09-02203-0}.

\bibitem{SloanWomersley2004Extremal}
{\sc I.~H. Sloan and R.~S. Womersley}, {\em Extremal systems of points and
  numerical integration on the sphere}, Adv. Comput. Math., 21 (2004),
  pp.~107--125, \url{https://doi.org/10.1023/B:ACOM.0000016428.25905.da}.

\bibitem{SloanWomersley2026QMC}
{\sc I.~H. Sloan and R.~S. Womersley}, {\em {QMC} designs---cubature on the
  sphere without polynomial exactness}, Constructive Approximation, 64 (2026),
  pp.~183--208, \url{https://doi.org/10.1007/s00365-026-09760-9}.

\bibitem{WuMartinsson2023}
{\sc B.~Wu and P.-G. Martinsson}, {\em A unified trapezoidal quadrature method
  for singular and hypersingular boundary integral operators on curved
  surfaces}, SIAM J. Numer. Anal., 61 (2023), pp.~2182--2208,
  \url{https://doi.org/10.1137/22M1520372}.

\bibitem{Wu2026PerturbedQMC}
{\sc H.-N. Wu}, {\em Perturbed {QMC} designs and reduced {QMC} strength}, arXiv
  preprint arXiv:2611.28831,  (2026).

\bibitem{WuZhuang2026NestedQMC}
{\sc H.-N. Wu and X.~Zhuang}, {\em Nested {QMC} designs on spheres}, arXiv
  preprint arXiv:2609.15767,  (2026).

\bibitem{WuZhangZhao2026}
{\sc Z.~Wu, Z.~Zhang, and X.~Zhao}, {\em Quasi-{Monte Carlo} time-splitting
  methods for the {Schr\"o}dinger equation with {Gaussian} random potential},
  J. Sci. Comput., 109 (2026), p.~85,
  \url{https://doi.org/10.1007/s10915-026-03461-y}.

\end{thebibliography}
\end{document}